\documentclass[12pt]{amsart}
\usepackage{amsfonts,amssymb,enumerate,fullpage,epsf}

\input xy
\xyoption{all}

\def\Pr{\begin{proof}}
\def\Rp{\end{proof}}

\def\IR{\mathbb R}
\def\IN{\mathbb N}

\def\ZF{\mathbf {ZF}}
\def\ZFC{\mathbf {ZFC}}

\def\DC{\mathbf {DC}}

\def\T{{\mathbf T}}

\def\AC{\mathbf {AC}}

\def\ACD{\mathbf {AC_{\IN}}}
\def\ACDF{\mathbf {AC_{\IN}^{fin}}}

\def\ZFC{\mathbf {ZFC}}

\def\AUc{\mathbf {A1}}
\def\AH{\mathbf {A2}}
\def\AHb{\mathbf {A3}}
\def\AHbf{\mathbf {A4}}

\def\Uwo{\mathbf {Uwo}}

 \DeclareMathOperator{\sgn}{sgn}

\theoremstyle{plain}

\newtheorem{corollary}{Corollary}
\newtheorem{proposition}{Proposition}
\newtheorem*{proposition*}{Proposition}
\newtheorem{theorem}{Theorem}
\newtheorem*{theorem*}{Theorem}
\newtheorem*{corollary*}{Corollary}
\newtheorem{lemma}{Lemma}
\newtheorem*{lemma*}{Lemma}

\theoremstyle{definition}
\newtheorem{definition}{Definition}
\newtheorem{notation}{Notation}
\newtheorem{question}{Question}

\theoremstyle{remark}
\newtheorem{remark}{Remark}
\newtheorem{example}{Example}

\date{\today}

\begin{document}
\title[Uniform Eberlein spaces]{Uniform Eberlein spaces and the Finite Axiom of Choice}
\author[M.~Morillon]{Marianne Morillon}
 \address{ERMIT, D\'epartement de Math\'ematiques et Informatique,
 Universit\'e de La R\'eunion, 15 avenue Ren\'e Cassin - BP 7151 -
 97715 Saint-Denis Messag. Cedex 9 FRANCE}
 \email[Marianne Morillon]{mar@univ-reunion.fr}
 \urladdr{http://personnel.univ-reunion.fr/mar}
 \subjclass[2000]{Primary 03E25~;  Secondary 54B10, 54D30, 46B26}
 \keywords{Axiom of Choice, product topology, compactness, Eberlein spaces, uniform Eberlein spaces}

\begin{abstract} 
We work in set-theory without choice $\ZF$. Given a closed subset $F$ of $[0,1]^I$ which is 
a bounded subset of $\ell^1(I)$ ({\em resp.}  such that $F \subseteq \ell^0(I)$),
we show that  the countable  axiom of choice for finite subsets of $I$, ({\em resp.} the countable axiom of choice $\ACD$)
  implies that $F$  is compact. 
This enhances previous results where $\ACD$ ({\em resp.} the axiom of Dependent Choices $\DC$) 
was required. Moreover, if $I$ is linearly orderable (for example $I=\IR$), the closed unit ball of $\ell^2(I)$ is
 weakly compact (in $\ZF$).
\end{abstract}

 \maketitle

 \begin{center}
    \mbox{\epsfbox{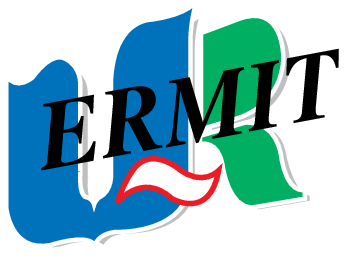}}
    \end{center}

 \begin{center}
 EQUIPE R\'EUNIONNAISE DE MATH\'EMATIQUES ET INFORMATIQUE TH\'EORIQUE (ERMIT)
 \end{center}
 \bigskip

 \bigskip

\tableofcontents

\section{Introduction} \label{sec:intro}
We work in the set-theory without the Axiom of Choice $\ZF$.
It is a well known theorem of Kelley (see \cite{Kel}) that, in $\ZF$,  the Axiom of Choice (for short $\AC$) is equivalent to the Tychonov axiom
$\T$: {\em ``Every family $(X_i)_{i \in I}$ of compact topological spaces  has a compact product.''}
Here, a topological space $X$ is {\em compact} if every family $(F_i)_{i \in I}$ of closed subsets of $X$ satisfying the finite intersection property ({\em FIP}) has a non-empty intersection. However, some particular cases of the Tychonov axiom are provable in $\ZF$, for example: 
\begin{remark} \label{rem:finite-prod-comp} A finite product of compact spaces is compact (in $\ZF$).
\end{remark}

 Say that the  topological space $X$ is {\em closely-compact} if there is a mapping $\Phi$ associating to every  family $(F_i)_{i \in I}$ of closed subsets of $X$ satisfying the 
{\em FIP} an element $\Phi((F_i)_{i \in I})$ of $\cap_{i \in I}F_i$: the mapping $\Phi$ is a {\em witness of closed-compactness} on $X$.
Notice that a compact topological space $X$ is closely-compact if and only if  there exists a mapping $\Psi$ associating to every non-empty closed subset $F$ of $X$ an element of $F$.

\begin{example} \label{ex:otc-comp} Given a linear order $(X,\le)$  which is complete (every non-empty subset of $X$ has a least upper bound), then the order topology on $X$ is closely compact. In particular, the closed bounded interval $[0,1]$
of $\IR$ is closely compact.
\end{example}
\Pr The space $X$ is compact (the classical proof is valid in $\ZF$). Moreover $X$ is closely compact since one can consider the choice function associating to every non-empty closed subset its first element. 
\Rp

The following Theorem is provable in $\ZF$:
\begin{theorem*}[\cite{Fo-Mo}] \label{theo*:prod-clos-comp}
Let $\alpha$ be an ordinal. 
 If $(X_i,\Phi_i)_{i \in \alpha}$ is a family of witnessed closely-compact spaces, then 
$\prod_{i \in \alpha} X_i$ is closely-compact, and has  a witness of closed-compactness
which is definable from $(X_i,\Phi_i)_{i \in \alpha}$.
\end{theorem*}

\begin{example} For every ordinal $\alpha$, the product topological space $[0,1]^{\alpha}$ is closely compact in $\ZF$.
\end{example}

Given a set $I$, denote   by $B_1(I)$ the set of $x=(x_i)_{i \in I} \in \IR^I$ such that $\sum_{i \in I} |x_i| \le 1$: then
$B_1(I)$ is a closed subset of $[-1,1]^I$. 
In this paper, we shall prove that $B_1(I)$ is compact  using the {\em  countable axiom of choice for finite subsets of $I$} (see Theorem~\ref{theo:ac-fin2-ball-hilb-base} in Section~\ref{subsec:acfin2-u-eber}). This enhances Corollary~1 of \cite{Mo07}
and   partially solves Question~2 in \cite{Mo07}. 
We shall  deduce  (see Corollary~\ref{cor:B1deLO}) that, if $I$ is  linearly orderable, 
every closed subset of $[0,1]^I$ which is contained in 
$B_1(I)$ is closely compact. In particular, the closed unit ball of the Hilbert space $\ell^2(\IR)$ is compact in $\ZF$, and
this solves Question~3 of \cite{Mo07}. Notice that $\{0,1\}^{\IR}$ (and $[0,1]^{\IR}$) is not compact in $\ZF$ (see \cite{Ker01}).
We shall also prove that  {\em Eberlein} closed subsets of $[0,1]^I$ are compact using the {\em countable axiom of choice for
subsets of $I$} (see Corollary~\ref{cor:acd2comp-eberlein}) of Section~\ref{subsec:acd2eberl}. This enhances Corollary~3  in \cite{Mo07} where the same result  was proved using the axiom of {\em Dependent Choices} $\DC$. This also solves Questions~4 and~5 thereof.

The paper is organized as follows: in Section~\ref{sec:various-consAC} we review various consequences of $\AC$ (in particular 
the {\em countable  axiom of choice $\ACD$} and  the axiom of choice {\em restricted to finite subsets $\ACDF$}) and the known links between them. In Section~\ref{sec:var-eberlein} we present definitions of uniform Eberlein spaces, 
strong Eberlein spaces and Eberlein spaces. In Section~\ref{sec:comp-ZF} we give some tools for compactness or
sequential compactness in $\ZF$. In Section~\ref{sec:alex-comp}, we recall the one-point compactification $\hat X$ of a discrete
space $X$, and we show that for every  ordinal $\alpha \ge 1$, the closed-compactness of ${\hat X}^{\alpha}$ is equivalent
to the axiom of choice restricted to finite subsets of $X$. Finally, in Section~\ref{subsec:ac-fin2compB1} ({\em resp.} \ref{sec:dc-and-comp})  we prove that the  countable axiom of choice for finite sets ({\em resp.} the countable axiom of choice) implies  that uniform Eberlein spaces ({\em resp.} Eberlein spaces) are closely compact ({\em resp.} compact.) A basic tool for these two last Sections is a ``dyadic representation'' of elements of powers of $[0,1]$ (see the Theorem in Section~\ref{subsec:dyadic}) which
we  found in \cite[Lemma~1.1]{B-R-W}, and for which the authors cite \cite{Sim76}.

\section{Some weak forms of $\AC$} \label{sec:various-consAC}
In this Section, we  review some weak forms of the Axiom of Choice  which  will be used in this paper and some  known links between them. 
For detailed references and much  information on this subject, see \cite{Ho-Ru}.

\subsection{Restricted axioms of choice}
 Given  a formula $\phi$ of set-theory with one free variable $x$, consider the following  consequence of $\AC$, denoted by  
$\AC(\phi)$: {\em ``For every non-empty  family $A=(A_i)_{i \in I}$ of non-empty sets such that  $\phi[x/A]$ holds, 
then $\prod_{i \in I}A(i)$ is non-empty.''}

\begin{notation} In the particular case where  the formula $\phi$ says that ``$x$ is a mapping with domain  $I$  
with values in some $\ZF$-definable class $\mathcal C$'', the statement $\AC(\phi)$ is denoted by $\AC_I^{\mathcal C}$.
\end{notation}

The statement $\forall I \AC_I^{\mathcal C}$ is denoted by $\AC^{\mathcal C}$.
The statement $\AC_I^{\mathcal C}$ where $\mathcal C$ is the collection of all sets is denoted by $\AC_I$.

\begin{notation} For every set $X$, we denote by $fin(X)$ the set of finite subsets of $X$. We denote by $fin$ the (definable) class of finite sets. 
\end{notation}
 
So, given a set $X$,   $\AC^{fin(X)}$ is the  following statement: 
 {\em ``For every non-empty  family $(F_i)_{i \in I}$ of non-empty finite subsets of $X$, $\prod_{i \in I} F_i$ is non-empty.''},  
and  $\AC^{fin}$ is the following statement: 
 {\em ``For every  non-empty family $(F_i)_{i \in I}$ of non-empty finite sets, $\prod_{i \in I} F_i$ is non-empty.''} 
The {\em countable Axiom of Choice} says that:
\begin{quote}$\ACD$: {\em If $(A_n)_{n \in \IN}$ is a family of  non-empty sets, then there exists a mapping $f : \IN \to \cup_{n \in \IN}A_n$ 
associating to every $n \in \IN$ an element $f(n) \in A_n$.}
\end{quote} 
And the {\em countable Axiom of Choice for finite sets} says that:
\begin{quote}$\ACDF$: {\em If $(A_n)_{n \in \IN}$ is a family of finite non-empty sets, then there exists a mapping $f : \IN \to \cup_{n \in \IN}A_n$ associating to every $n \in \IN$ an element $f(n) \in A_n$.}
\end{quote}

\subsection{Well-orderable union of finite sets}
Given an infinite ordinal $\alpha$, and a class $\mathcal C$ of sets, we consider  the following consequence of  $\AC^{\mathcal C}$:
\begin{quote} $\Uwo_{\alpha}^{\mathcal C}$: {\em For every family $(F_i)_{i \in \alpha}$ of elements of $\mathcal C$, 
the set $\cup_{i \in \alpha} F_i$ is well-orderable.}
\end{quote}

\begin{remark}
 $\AC^{fin(X)}$ implies $\Uwo_{\alpha}^{fin(X)}$. 
\end{remark}

\subsection{Dependent Choices} The axiom of {\em Dependent Choices} says  that:
\begin{quote} $\DC$: {\em Given a non-empty set $X$ and a binary relation $R$ on $X$ such that $\forall x \in X \exists y \in X \; xRy$, then there exists a sequence 
$(x_n)_{n \in \IN}$ of $X$ such that for every $n \in \IN$, $x_n R x_{n+1}$.}
\end{quote} 

Of course, $\AC \Rightarrow \DC \Rightarrow \ACD \Rightarrow \ACDF$. However, the converse statements are not provable in $\ZF$, and $\ACDF$ is not provable in $\ZF$ (see references in \cite{Ho-Ru}).

\subsection{The ``Tychonov'' axiom} \label{subsec:tycho}
Given a class $\mathcal C$ of  compact topological spaces and a set $I$, we consider the following consequence of the Tychonov axiom:
\begin{quote} $\T_I^{\mathcal C}$: {\em Every family $(X_i)_{i \in I}$ of spaces belonging to the class $\mathcal C$ has a compact product.} 
\end{quote} 

For example $\T_{\IN}^{fin(X)}$ is the statement {\em ``Every sequence of finite discrete subsets of $X$ has a compact product.''}

\begin{remark} \label{rem:acfin-alpha} 
\begin{enumerate}[(i)]
\item \label{it:rem-acfin1} Given a set $X$, for every ordinal $\alpha$, 
$$\AC^{fin(X)} \Rightarrow \Uwo_{\alpha}^{fin(X)} \Rightarrow \T_{\alpha}^{fin(X)} \Rightarrow \AC_{\alpha}^{fin(X)}$$
\item \label{it:rem-acfin2} For every ordinal $\alpha$, $\Uwo_{\alpha}^{fin} \Leftrightarrow  \T_{\alpha}^{fin} \Leftrightarrow \AC_{\alpha}^{fin}$.
\end{enumerate}
\end{remark}
\Pr \eqref{it:rem-acfin1} $\Uwo_{\alpha}^{fin(X)} \Rightarrow \T_{\alpha}^{fin(X)}$: Given a family $(F_i)_{i \in \alpha}$ of finite subsets of $X$, the
statement $\Uwo_{\alpha}^{fin(X)}$ implies the existence of a family $(\Phi_i)_{i \in \alpha}$ such that for each $i \in \alpha$, the
discrete space is closely compact with witness $\Phi_i$. Using the  Theorem of Section~\ref{sec:intro}, it follows that $\prod_{i \in \alpha}F_i$
is (closely) compact. 
$ \T_{\alpha}^{fin(X)} \Rightarrow \AC_{\alpha}^{fin(X)}$: one can use Kelley's argument (see \cite{Kel}). \\
\eqref{it:rem-acfin2} For $\AC_{\alpha}^{fin} \Rightarrow \Uwo_{\alpha}^{fin}$: given some family $(F_i)_{i \in \alpha}$ of
finite non-empty sets, then, for each $i \in \alpha$, denote by $c_i:=\{0..c_{i-1}\}$ the (finite) cardinal of $F_i$; 
thus set $G_i$ of one-to-one mappings from $F_i$ to $c_i$ is finitel, and, by  $\AC_{\alpha}^{fin}$, the set 
$\prod_{i \in \alpha}G_i$ is non-empty. This implies a well-order on the set $\cup_{i \in \alpha}F_i$. 
\Rp

\section{Some classes of closed subsets of $[0,1]^I$} \label{sec:var-eberlein}
\begin{notation}
Let $I$ be a set. Given  some element $x=(x_i)_{i \in I} \in \IR^I$,
denote by $supp(x)$ the {\em support} $\{i \in I : x_i \neq 0\}$. 
Given some subset $A$ of $\IR$ containing $0$, denote by $A^{(I)}$ the set of elements of $A^I$ with {\em finite} support. 
\end{notation} 

We endow the  space $\IR^I$ with the  product topology, which we denote
by  $\mathcal T_I$.

\subsection{Eberlein closed subsets of $[0,1]^I$}
Given a set $I$, we denote by $\ell^{\infty}(I)$ the Banach space of bounded mappings $f:I \to \IR$, endowed with the ``sup'' norm. If $I$ is infinite,
we denote
by $c_0(I)$ the closed subspace of $\ell^{\infty}(I)$ consisting of $f \in \ell^{\infty}(I)$ such that $f$ converges to $0$ according to the Fr\'echet
filter on $I$ ({\em i.e.} the set of cofinite subsets of $I$). Thus
$$\ell^0(I):=\{x=(x_i)_{i \in I} : \; \forall \varepsilon>0 \exists F_0 \in \mathcal P_f(I) \forall i \in I \backslash F_0 \;
|x_i| \le \varepsilon\}$$ 
If $I$ is finite, then we define $c_0(I):=\ell^{\infty}(I)=\IR^I$. 

\begin{definition} A  topological space $F$ is {\em $I$-Eberlein} if $F$ is a closed subset of $[0,1]^I$ and if   
$F \subseteq c_0(I)$. A topological space $X$ is {\em Eberlein} if $X$ is  homeomorphic with some $I$-Eberlein space.  
\end{definition}

\begin{remark}
 Amir and Lindenstrauss (\cite{Ami-Lin}) proved in $\ZFC$ that  every  weakly compact subset of a normed space is an Eberlein space. 
This result relies on the existence of a Markhushevich basis in every weakly compactly generated Banach space, and the proof of the existence of
such a basis (see \cite{Fab-et-al}) relies on (much)  Axiom of Choice. 
\end{remark}

\begin{remark} Consider the compact topological space $X:=[0,1]^{\IN}$. Then, the closed subset $X$ of $[0,1]^{\IN}$
is not $\IN$-Eberlein. However, the mapping 
$f : X \to [0,1]^{\IN} \cap c_0(\IN)$
associating to each $x=(x_n)_{n \in \IN} \in X$ the element $(\frac{x_n}{n+1})_{n \in \IN}$ is  continuous and one-to-one, so
$X$ is homeomorphic with the compact (hence closed) subset $f[X]$ of $[0,1]^{\IN} \cap c_0(\IN)$. It follows that 
$X$ is homeomorphic with some $\IN$-Eberlein space. 
\end{remark}

\begin{proposition} \label{prop:Eberlein-ppties}
\begin{enumerate}[(i)]
\item \label{it:Eber1} Every closed subset of a $I$-Eberlein ({\em resp.} Eberlein) space is $I$-Eberlein ({\em resp.} Eberlein).  
\item \label{it:Eber2}  Let $(I_n)_{n \in \IN}$ be a sequence of pairwise disjoint sets, and denote by $I$ the set $\sqcup_{n \in \IN}I_n$.
Let $(F_n)_{n \in \IN}$ be a sequence of topological spaces such that each $F_n$ is $I_n$-Eberlein . Then the closed subset 
$\prod_{n \in \IN} F_n$ of $[0,1]^I$ is homeomorphic with a $I$-Eberlein space. 
\end{enumerate}
\end{proposition}
\Pr \eqref{it:Eber1} is trivial. We prove  \eqref{it:Eber2}. For every $n \in \IN$, let $f_n : F_n \to [0,1]^{I_n}$ be the mapping associating to each
$x \in F_n$  the element $\frac{1}{n+1}f_n(x)$ of $[0,1]^{I_n}$. Let  $f:=\prod_{n \in \IN }f_n : \prod_{n \in \IN} F_n \to [0,1]^I$.
Then $f$ is one-to-one and continuous. Moreover, the  subset $F:=Im(f)$ of $[0,1]^I$ is closed since $F$ 
is the product $\prod_{n \in \IN} \tilde F_n$
where for each $n \in \IN$, $\tilde F_n$ is the closed subset $\frac{1}{n+1}.F_n$ of $[0,1]^{I_n}$. Finally, it can be easily checked
that  $F \subseteq c_0(I)$. 
\Rp

\subsection{Uniform Eberlein closed subsets of $[0,1]^I$}
\subsubsection{The ball $B_p(I)$, for $1 \le p < +\infty$.}
For every real number $p \ge 1$, define as usual the 
normed space  $\ell^p(I):=\{(x_i)_{i \in I} : \sum_i |x_i|^p < +\infty\}$ endowed with the  norm 
$N_p : x=(x_i)_{i \in I} \mapsto ({\sum_{i}|x_i|^p})^{1/p}$. We denote by $B_p(I)$ the large
unit ball $\{x \in \IR^I : \sum_i |x_i|^p \le 1\}$ of $\ell^p(I)$. 
Notice that for $p=1$ ({\em resp.} $1<p<+\infty$) the topology induced by $\mathcal T_I$ on   $B_p(I)$ is   the topology induced by the weak* topology $\sigma(\ell^1(I),\ell^0(I))$ ({\em resp.}  the topology
induced by the weak topology $\sigma(\ell^p(I),\ell^q(I))$ where $q=\frac{p}{p-1}$ is the conjuguate of $p$). Also notice that
for $1 \le p < +\infty$, $B_p(I)$ is a closed subset of $[0,1]^I$.

\begin{proposition}  \label{prop:B1-to-Bp} If  $1 \le p< +\infty$, then   $B_p(I)$ is homeomorphic with $B_1(I)$.
\end{proposition}
\Pr Consider the mapping  $h_p: B_1(I) \to B_p(I)$ associating to every $x=(x_i)_i \in B_1(I)$ the family  $(\sgn(x_i) |x_i|^{1/p})_{i \in I}$.
\Rp 

It follows that for every $p,q \in [1,+\infty[$, spaces $B_p(I)$ and $B_q(I)$ are homeomorphic {\em via} 
$h_{p,q}:=h_q \circ h_p^{-1}: B_p(I) \to B_q(I)$.

\subsubsection{Uniform Eberlein spaces}
Given a set $I$, and some real number  $p \in [1,+\infty[$, 
we denote by  $B^+_p(I)$ the positive ball of $\ell^p(I)$:
$$B^+_p(I) := \{x=(x_i)_{i \in I} \in [0,1]^I : \sum_{i \in I} x_i^p \le 1\}$$

\begin{definition} A  topological space $F$ is {\em $I$-uniform Eberlein} if 
there exists a real number $p \in [1,+\infty[$ such that $F$ is   a closed subset of   $B^+_p(I)$. 
 A topological space $X$ is {\em uniform Eberlein} if $X$ is  homeomorphic with some $I$-uniform Eberlein space.  
\end{definition}

Of course, every $I$-uniform Eberlein space is $I$-Eberlein. Moreover, using Proposition~\ref{prop:B1-to-Bp},
 every $I$-uniform Eberlein space is homeomorphic with a closed
subset of $B^+_1(I)$.

\begin{proposition} \label{prop:unif-eber-prod-den}
\begin{enumerate}[(i)] 
\item \label{it:unif-Eberl1} For every set $I$, every closed subset of a $I$-uniform Eberlein space  is $I$-uniform Eberlein.
\item  \label{it:unif-Eberl2}  Let $(I_n)_{n \in \IN}$ be a sequence of pairwise disjoint sets, and denote by $I$ the set $\sqcup_{n \in \IN}I_n$.
Let $(F_n)_{n \in \IN}$ be a sequence of topological spaces such that each $F_n$ is a $I_n$-uniform Eberlein space. Then the closed subset 
$F:=\prod_{n \in \IN} F_n$ of $[0,1]^I$ is $I$-uniform Eberlein. 
\end{enumerate}
\end{proposition}
\Pr \eqref{it:unif-Eberl1} is easy. The proof of \eqref{it:unif-Eberl2} is similar to the proof of 
Proposition~\ref{prop:Eberlein-ppties}-\eqref{it:Eber2}.
\Rp

In particular, the compact space $[0,1]^{\IN}$ (and thus every metrisable compact space) is $\IN$-uniform Eberlein.
For every set $I$,  ${B^+_1(I)}^{\IN}$ is $(I \times \IN)$-uniform Eberlein.

\begin{remark}
Let    $Z:=\cap_{i \in I} \{(x,y) \in B^+_1(I) \times B^+_1(I) : \; x_i .y_i=0\}$: then $Z$ is a closed subset of $B^+_1(I) \times B^+_1(I)$,  
and  the mapping 
$- : Z \to B_1(I)$ is an homeomorphism; it follows that $B_1(I)$ is homeomorphic with a  $(I \times  \{0,1\})$-uniform Eberlein space.  
\end{remark}

\subsubsection{Weakly closed bounded subsets  of a Hilbert space}
\begin{remark}
Given a Hilbert space $H$ with a Hilbert basis $(e_i)_{i \in I}$, then its closed unit ball
(and thus every bounded weakly closed subset of $H$) is (linearly) homeomorphic with
the uniform Eberlein space $B_2(I)$. 
\end{remark}

\medskip

Consider the  following statements (the first two ones were introduced in \cite{De-Mo} and \cite{Mo04} 
and are consequences of the Alaoglu theorem):  
\begin{itemize}
\item $ \AUc$: The closed unit ball (and thus every bounded subset which is closed in the convex topology) of a uniformly convex Banach space is compact in the convex topology.
\item $ \AH$: (Hilbert) The closed unit ball (and thus every bounded weakly closed subset) of a Hilbert space is weakly compact. 
\item $ \AHb$: (Hilbert with hilbertian basis) For every set $I$, the closed unit ball of $\ell^2(I)$ is 
weakly compact. 
\item $ \AHbf$:  For every sequence $(F_n)_{n \in \IN} $ of finite sets, the closed unit ball of 
$\ell^2(\cup_{n \in \IN} F_n)$ is weakly compact. 
\end{itemize}
Of course, $ \AUc \Rightarrow  \AH \Rightarrow  \AHb \Rightarrow  \AHbf$.

\begin{theorem*}[\cite{Fo-Mo}, \cite{Mo07}] \label{theo:acd2wc-uc} 
\begin{enumerate}[(i)] 
\item \label{it:ACD2A1} $\ACD \Rightarrow  \AUc$. 
\item \label{it:A1not2ACD}  $ \AUc \not \Rightarrow \ACD$.
\item \label{it:ACD2A2} $ \AHbf \Rightarrow \ACDF$.    
\end{enumerate}
\end{theorem*}

In this paper, we will prove that the following statements are equivalent: $\AHb$, $\AHbf$, 
$\AC^{fin}_{\omega}$  (see Corollary~\ref{cor:acdf&Ahb}). 

\begin{question}
Does   $\AH$ imply $\AUc$? Does $\AHb$ imply $\AH$? 
\end{question}

\begin{remark} \label{rem:ah-wo} If a Hilbert space $H$ has a well orderable dense subset, then $H$ has a well orderable hilbertian basis, thus $H$ is isometrically isomorphic with some
$\ell^2(\alpha)$ where $\alpha$ is an ordinal. In this case, the closed unit ball of $H$ endowed with the weak topology 
is homeomorphic with a closed subset of $[-1,1]^{\alpha}$, so  this ball 
 is weakly compact. 
\end{remark}

\subsection{Strongly Eberlein closed subsets of $[0,1]^I$}
\begin{definition} A  topological space $F$ is {\em $I$-strong  Eberlein} if $F$ is a closed subset of $[0,1]^I$
which is contained in $\{0,1\}^{(I)}$.  
 A topological space $X$ is {\em strong Eberlein} if $X$ is  homeomorphic with some $I$-strong Eberlein space.  
\end{definition}

Of course, every $I$-strong Eberlein set is $I$-Eberlein.

\begin{remark}
For every set $I$, every closed subset of a $I$-strong  Eberlein space  is $I$-strong  Eberlein.
\end{remark}

\section{Compactness (in $\ZF$)} \label{sec:comp-ZF}
\subsection{Lattices and filters} \label{subsec:filters}
 Given a lattice $\mathcal L$ of subsets of a set $X$, say that  a non-empty proper subset $\mathcal F$  of $\mathcal L$ is a {\em filter} if  it satisfies the two following conditions:
\begin{enumerate}[(i)]
 \item $\forall A, B \in \mathcal F, \; A \cap B \in \mathcal F$
\item $\forall A \in \mathcal F, \; \forall B \in \mathcal L, (A \subseteq B \Rightarrow B \in \mathcal F)$
\end{enumerate}
Say that an element  $A \in \mathcal L$ is {\em $\mathcal F$-stationar} if for every $F \in \mathcal F$, $A \cap F \neq \varnothing$.

\begin{remark} \label{rem:stat-comp}
 Let $X$ be a topological space, let  $\mathcal L$ be a lattice of closed subsets of $X$,
and let $\mathcal F$ be a filter of $\mathcal L$. Let $K \in \mathcal L$. 
If $K$ is a compact subset of $X$ and if $K$ is $\mathcal F$-stationar, then  $\cap \mathcal F$ is non-empty. 
\end{remark}

\begin{definition}
Given a family $(X_i)_{i \in I}$ of topological spaces, and denoting by $X$ the topological product of this family, a closed subset $F$ of $X$ is {\em elementary} if $F$ is a  finite union of sets of
the form $\prod_{i \neq i_0}X_i \times C$ where $i_0 \in I$ and 
$C$ is a closed subset of $X_{i_0}$.  
\end{definition}

Given a family $(X_i)_{i \in I}$ of topological spaces with product $X$,
the set of elementary closed subsets of $X$ is a lattice of subsets of $X$ that we denote by $\mathcal L_X$. Notice that given a elementary closed subset $F$ of $X$, and some subset $J$ of $I$, the projection $p_J[F]$ is a closed subset of $\prod_{j \in J}X_j$.

\subsection{Continuous image of a compact space} \label{subsec:comp&close-comp}
The following Proposition is easy:
\begin{proposition} \label{prop-im-cont}
 Let $X$, $Y$ be topological spaces and let $f: X \twoheadrightarrow Y$ be a continuous onto mapping. If $X$ is compact ({\em resp.} closely-compact),
 then $Y$ is also compact ({\em resp.} closely compact). If $\Phi$ is a witness of closed-compactness on $X$, 
then $Y$ is closely-compact, and has a witness of closed-compactness which is definable from
$f$ and $\Phi$. 
\end{proposition}

\subsection{Sequential compactness} \label{subsec:seq-compact}
We denote by $[\IN]^{\omega}$ the set of infinite subsets of $\IN$.

\begin{definition}
 A topological space $X$ is {\em sequentially compact} if every sequence $(x_n)_{n \in \IN}$ of $X$
has an infinite subsequence which converges in $X$. A {\em witness of sequential compactness}
on $X$ is a mapping $\Phi : X^{\IN} \to [\omega]^{\omega} \times X $ associating to each sequence $(x_n)_{n \in \IN}$ of $X$ an element $(A,l) \in [\omega]^{\omega} \times X$ such
that  $(x_n)_{n \in A}$ converges to $l$.
\end{definition}

\begin{example} If $(X,\le)$ is a complete linear  order, then $X$ is sequentially compact, with a witness definable from
$(X,\le)$: given a sequence $(x_n)_{n \in \IN}$, build some infinite subset $A$ of $\IN$ such that $(x_n)_{n \in A}$ is monotone;
then if  $(x_n)_{n \in A}$ is ascending ({\em resp.} descending), then $(x_n)_{n \in A}$ converges to $\sup_{n \in A} x_n$
({\em resp.} $\inf_{n \in A} x_n$).
\end{example}

\begin{example} \label{ex:top-kelley} Given an infinite set $X$, and some set $\infty \notin X$, consider the topology on 
$\tilde X:=X \cup \{\infty\}$
generated by cofinite subsets of $\tilde X$ and $\{\infty\}$. This topology is compact and $T_1$ but it is not $T_2$. 
This topology is sequentially compact, and, given a point $a \in X$, there is a  witness of sequential compactness which is definable from $X$ and $a$: given a sequence $(x_n)_{n \in \IN}$ of $\tilde X$, either the set of terms 
$\{x_n : n \in \IN\}$ is finite, and then one can define by induction an infinite subset $A$ of $\IN$ such that 
$\{x_n : n \in A\}$ is constant; else one can  define by induction an infinite subset $A$ of $\IN$ such that 
$\{x_n : n \in A\}$ is one-to-one, thus it converges to $a$ (and also to every point in $X$).
\end{example}

Notice that the topology in Example~\ref{ex:top-kelley} is the one used by Kelley (see \cite{Kel}) to prove that ``Tychonov implies $\AC$''. The following Lemma is easy:
\begin{lemma} \label{lem:imag-cont-seq-comp}
 Let $X,Y$ be two topological spaces and let $f: X  \twoheadrightarrow Y$ be an onto continuous mapping
which has a section $j$ (for example if $f$ is one-to-one).
If $X$ is sequentially compact, then  $Y$ is also sequentially compact. Moreover, if
there is a witness $\Phi$ of sequential compactness on $X$, there also exists a witness of sequential
compactness on $Y$ which is definable from $f$,$\Phi$ and $j$. 
\end{lemma}

\begin{lemma} \label{lem:prod-wit-seq-comp}
 Let $(X_n,\phi_n)_{n \in \IN}$ be a sequence of witnessed sequentially compact spaces.
 The space $\prod_{n \in \IN} X_n$ is sequentially compact, and has a witness definable from 
$(X_n,\phi_n)_{n \in \IN}$. 
\end{lemma}
\Pr Usual diagonalization.
\Rp 

\begin{example} \label{ex:seq-compact} If $D$ is a countable set, then the topological space $[0,1]^D$ is sequentially compact,
a witness of sequential compactness beeing definable from every well order on $D$.
\end{example}

Say that  a sequentially compact topological space $X$ is {\em witnessable} if there exists a witness of sequential compactness
on $X$. 
It follows from Lemma~\ref{lem:prod-wit-seq-comp}, that with $\ACD$, every sequence   $(K_n)_{n \in \IN}$ of  
 witnessable sequentially compact spaces has a product which is sequentially compact.

\subsection{$\ACD$ and countable products of compact spaces}
Denote by $T^{comp}_{\omega}$ the following statement: {\em ``Every sequence of compact spaces has a compact product.''}
Then Kelley's argument shows that $T^{comp}_{\omega} \Rightarrow \ACD$. However, 
 it is an open question (see \cite{Br85}, \cite{HKRS}) to know whether $\ACD$ implies $T^{comp}_{\omega}$.

\begin{definition}
 A topological space $X$ is {\em $\omega$-compact} if every descending sequence 
$(F_n)_{n \in \IN}$ of non-empty closed subsets of $X$ has a non-empty intersection.
Say that the space $X$ is {\em cluster-compact} if every sequence $(x_n)_{ \in \IN}$ of $X$ has a {\em cluster point} {\em i.e.} the set
$\cap_{n \in \IN} \overline{\{x_k : k \ge n\}}$ is non-empty.
\end{definition}

\begin{remark} \label{rem:acd-seq-comp}
\begin{enumerate}[(i)]
\item \label{it:acd-seq-comp1} Notice that $\text{sequentially compact}  \Rightarrow \text{``cluster-compact''}$. Also notice that  
 $\text{``$\omega$-compact''}  \Rightarrow \text{``cluster-compact''}$ and that the converse holds with $\ACD$ 
(see \cite[Lemma~1]{HKRS}). 
\item Given a sequence  $(K_n)_{n \in \IN}$ of compact spaces, then, denoting by $K$ the product of this family,
 $K$ is compact {\em iff} $K$ is $\omega$-compact (see \cite[Theorem~6]{HKRS}). 
\item \label{it:acd-prod-seq-comp} If the product  $K$ of a sequence  $(K_n)_{n \in \IN}$ of compact spaces is
sequentially compact, then $\ACD$ implies that $K$ is compact. 
\end{enumerate}
\end{remark}

\begin{proposition} \label{prop:acd-prod-comp} $\ACD$ is equivalent to the following statement:
 {\em ``Every sequence $(K_n)_{n \in \IN}$ of witnessable sequentially compact spaces which are also compact  has a compact product.''}
\end{proposition}
\Pr $\Rightarrow$:  Given a sequence $(K_n)_{n \in \IN}$ of witnessable sequentially compact spaces which are also compact, then, 
using $\ACD$, one can choose a witness of sequential compactness on every space $K_n$. 
It follows by Lemma~\ref{lem:prod-wit-seq-comp} that  $K$ is sequentially compact, 
whence $K$ is compact by Remark~\ref{rem:acd-seq-comp}-\eqref{it:acd-prod-seq-comp}.   \\
$\Leftarrow$:  We use Kelley's argument (see \cite{Kel}). Let $(A_n)_{n \in \IN}$ be a sequence of non-empty sets. 
Consider some element $\infty \notin \cup_{n \in \IN} A_n$, and for every $n \in \IN$, 
 denote by $K_n$ the set $A_n \cup \{\infty\}$ endowed with the topology generated by $\{\infty\}$ and cofinite subsets of $K_n$
(see Example~\ref{ex:top-kelley}). Then each $K_n$ is compact and sequentially compact; moreover, given an element $a \in A_n$, there is a witness of sequential compactness on $K_n$ which is definable
from $A_n$,$\infty$ and $a$. So each $K_n$ is a witnessable sequentially compact space. 
It follows  from the hypothesis  that the product $K:=\prod_{n \in \IN}K_n$ is  compact. We end as in Kelley's proof: for every $n \in \IN$, let $F_n$ be the closed set 
$A_n \times \prod_{i \neq n} K_i$. By compactness of $K$, the set $\cap_{n \in \IN} F_n$ is non-empty. 
This yields an element of $\prod_{n \in \IN} A_n$.
\Rp

\section{One-point compactifications and related spaces} \label{sec:alex-comp}
\subsection{The one-point compactification of a set} \label{subsec:alex}
Given a set $X$,  we denote by  $\hat X$ the Alexandrov compactification of the (Hausdorff locally compact) discrete space  $X$: 
$\hat X := X \cup \{\infty\}$ where $\infty$ is some set $\notin X$ (for example 
$\infty:=\{x \in X : x \notin x\}$; if $X$ is finite, then $\hat X$ is discrete else open subsets of the space $\hat X$ are subsets of $X$ or cofinite subsets of $\hat X$ containing $\infty$. Notice that the  space $\hat X$ is compact and Hausdorff in $\ZF$.

\begin{example} \label{ex:hatX_unif-eb} Given a discrete topological space $X$, the one-point compactification ${\hat{X}}$ of $X$ is 
$X$-uniform Eberlein: consider  the Hilbert space $\ell^2(X)$; and denote by $(e_i)_{i \in X}$ the canonical basis of the vector space  
$\IR^{(X)}$;  then the subspace  $X=\{e_i : i \in X\}$ of $\IR^{(X)}$  is discrete and 
  the weakly closed and bounded subset $X \cup 0_{\IR^X}$ is the one-point compactification $\hat X$ of   $X$.
 \end{example}

\subsection{Various notions of compactness for $\hat X^{\alpha}$, $\alpha$ ordinal}
\subsubsection{$\hat X^{\IN}$ is  sequentially compact}
\begin{proposition} \label{prop:hatX-seq-comp} Let $X$ be an infinite set. 
\begin{enumerate}[(i)]
\item \label{it:X-chap-2} The space  $\hat X$ is  sequentially compact and has a witness of sequential compactness, definable from $X$. 
\item \label{polyad-den-seq} The space $\hat X^{\IN}$ is sequentially compact with a witness definable from $X$.
\end{enumerate}
\end{proposition}
\Pr \eqref{it:X-chap-2} We define a witness $\Phi$ of sequential 
compactness on $X$ as follows: given a sequence $x=(x_n)_{n \in \IN}$ of $\hat X$, if the set
$T:=\{x_k : k \in \IN\}$ is infinite, we build (by induction) some infinite subset $A$ of $\IN$
such that $\{x_k : k \in A\}$ is one-to-one, and we define  $\Phi(x):=(A,\infty)$;  else the set  
$T$ is finite, so    we build by induction some
infinite subset $A$ of $\IN$ such that the sequence $\{x_k : k \in A\}$ is a singleton $\{l\}$, and we define  $\Phi(x):=(A,l)$. \\
 \eqref{polyad-den-seq} We apply \eqref{it:X-chap-2} and Lemma~\ref{lem:prod-wit-seq-comp} in Section~\ref{subsec:seq-compact}.
\Rp

\subsubsection{$\AC^{fin}$ and closed-compactness}
\begin{proposition} \label{prop:X-hat-comp} Let $X$ be a set. 
\begin{enumerate}[(i)]
\item \label{it:X-hat-0} There is a mapping associating to every non-empty closed subset $F$ of $\hat X$, 
a finite non-empty closed subset $\tilde F$ of $F$. 
\item \label{it:X-hat-2} $\AC^{fin(X)} \Leftrightarrow \text{ The space $\hat X$ is  closely compact}$.  
\end{enumerate}
\end{proposition}
\Pr We may assume that $X$ is infinite. \\
 \eqref{it:X-hat-0} Given a non-empty closed subset $F$ of  $\hat X$, define  $\tilde F:= \{\infty\}$ if $\infty \in F$ and 
$\tilde F:=F$ if $F$ is finite and $\infty \notin F$. \\
\eqref{it:X-hat-2} Use \eqref{it:X-hat-0}. 
\Rp

\subsection{Spaces  ${\hat X}^{\alpha}$, $\alpha$ ordinal}
\begin{remark}
For every set $X$, the space $\hat X$ is $X$-uniform Eberlein, so, given an ordinal $\alpha$, 
${\hat X}^{\alpha}$ is $X \times \alpha$-uniform Eberlein  
(see Proposition~\ref{prop:unif-eber-prod-den}-\eqref{it:unif-Eberl2}).
\end{remark}

\begin{proposition} \label{theo:acdf-comp-polyad}
 Let $X$ be a set. Let $\alpha$ be an ordinal $\ge 1$. 
\begin{enumerate}[(i)]
\item \label{it:hatx-puiss-alpha-1}  $\T_{\alpha}^{fin(X)} \Leftrightarrow \text{``$\hat X^{\alpha}$  is  compact''}$. 
\item \label{it:hatx-puiss-alpha-2} $\AC^{fin(X)} \Leftrightarrow \text{``$\hat X^{\alpha}$  is  closely compact''}$. 
\end{enumerate}
\end{proposition}
\Pr \eqref{it:hatx-puiss-alpha-1} $\Rightarrow$:   Let  $P$ be the topological product space  $\hat X^{\alpha}$. Let  $\mathcal F$ be a filter of the lattice $\mathcal L_X$ of elementary closed subsets of  $P$. We are going to define  by transfinite recursion  a family   
$(G_n)_{n \in \alpha}$ of finite subsets of  $\hat X$ such that, denoting for every  $n \in \alpha$ by $Z_n$ the elementary closed subset 
$G_n \times {\hat X}^{\alpha \backslash \{n\}}$ of $P$, the set 
 $\mathcal F \cup \{Z_i : i < n\}$ satisfies the  finite intersection property.
Given some $n \in \alpha$,  we define $G_n$  in function of  $(G_i)_{i < n}$ as follows:  
denote by $\mathcal G$  the filter generated by $\mathcal F \cup \{Z_i : i <n\}$; since  $\hat X$ is compact,
the closed subset  $F_n:=\cap p_{\{n\}}[\mathcal G]$ is non-empty, so let 
 $G_n:=\tilde F_n$ and   let $Z_n:=G_n \times {\hat X}^{\alpha \backslash \{n\}}$.
 Denote by  $\tilde {\mathcal F}$ the filter  generated by $\mathcal F \cup \{Z_i : i <\alpha\}$. 
Using $\T_{\alpha}^{fin(X)}$, the product space $F:=\prod_{n \in \alpha} G_n$ is compact, and non-empty since 
$\T_{\alpha}^{fin(X)}$ implies $\T_{\alpha}^{fin(X)}$ (see Remark~\ref{rem:acfin-alpha}). 
Moreover, the closed subset $F$  of $P$ is  $\tilde{\mathcal F}$-stationnar: 
it follows from Remark~\ref{rem:stat-comp} of Section~\ref{subsec:filters} 
 that $\cap \tilde {\mathcal F}$ is non-empty, whence  $\cap \mathcal F \neq \varnothing$. 
$\Leftarrow$: Let $(F_i)_{i < \alpha}$ be a family of finite subsets of $X$, endowed with the discrete topology. Then $\prod_{i < \alpha} F_i$
is compact because it is a closed subset of the compact Hausdorff space ${\hat X}^{\alpha}$.  \\
 \eqref{it:hatx-puiss-alpha-2} $\Rightarrow$: Use Proposition~\ref{prop:X-hat-comp}-\eqref{it:X-hat-2} and the Theorem of Section~\ref{subsec:comp&close-comp}. 
$\Leftarrow$: If $\hat X^{\alpha}$  is  closely compact, then so is its continuous image $\hat X$, whence $\AC^{fin(X)}$ holds (using
Proposition~\ref{prop:X-hat-comp}-\eqref{it:X-hat-2}). 
\Rp

\subsection{Spaces $\sigma_n(X)$, $n$ integer $\ge 1$}
\begin{notation} Given a set $X$, for  every integer  $n \ge 1$, let 
$$\sigma_n(X):=\{x \in \{0,1\}^{(X)}: |supp(x)| \le n\}$$
\end{notation}
Thus $\sigma_n(X)$ is  the set of  elements of $\IR^{(X)}$ with support having at most $n$ elements.  
Notice that the space $\sigma_n(X)$ is strong Eberlein.

\begin{remark} \label{rem:sigma-n}
\begin{enumerate}[(i)]
\item \label{it:sigma1-ue} The space  $\sigma_1(X)$ is the one-point compactification of  the discrete space  $X$ (thus $\sigma_1(X)$ is uniform Eberlein).
\item \label{it:sigma-n-im-sigma1} The mapping 
$\mathcal U_n: (\sigma_1(X))^n \to \sigma_n(X)$ associating to each 
$(x_1, \dots,x_n)$ the set  $\cup_{1 \le i \le n} x_i$ is continuous. 
\end{enumerate}
\end{remark}
\Pr  \eqref{it:sigma1-ue} Use Example~\ref{ex:hatX_unif-eb}. \eqref{it:sigma-n-im-sigma1}: easy.
\Rp 

\subsubsection{Compactness and closed compactness of $\sigma_n(X)$}
\begin{proposition} \label{prop:sigma-n-comp}
 Let $X$ be a set, and let $n$ be some integer $ \ge 1$. 
\begin{enumerate}[(i)]
\item \label{it:sigma-n-comp1} Both  spaces $(\sigma_1(X))^n$ and $\sigma_n(X)$ are  compact.
\item \label{it:sigma-n-comp2} With $\AC^{fin(X)}$, both  spaces $(\sigma_1(X))^n$ and $\sigma_n(X)$ are closely compact 
(with witnesses of closed compactness
definable from $X$,  $n$ and some choice function on non-empty finite subsets of $X$).
\end{enumerate}
 \end{proposition}
\Pr The results \eqref{it:sigma-n-comp1} and \eqref{it:sigma-n-comp2} for $\sigma_n(X)$ follow from  the result  on $(\sigma_1(X))^n$  thanks to 
Proposition~\ref{prop-im-cont} and the continuous onto mapping 
$\mathcal U_n: (\sigma_1(X))^n \to \sigma_n(X)$ defined in Remark~\ref{rem:sigma-n}. The result  \eqref{it:sigma-n-comp1}
for $(\sigma_1(X))^n$ comes from Remark~\ref{rem:finite-prod-comp}. We prove \eqref{it:sigma-n-comp2}
for $(\sigma_1(X))^n$: with  $\AC^{fin(X)}$,  $\sigma_1(X)$ is closely compact, so 
the space  $(\sigma_1(X))^n$  is also closely compact because 
it is a finite power of a   closely compact space (use the Theorem of Section~\ref{sec:intro}).    
\Rp

\subsubsection{Sequential compactness of $\prod_{n \in \IN}\sigma_n(X)$}
\begin{proposition} \label{prop:seq-comp}
 Let $X$ be a set, and let $n$ be some integer $ \ge 1$. 
\begin{enumerate}[(i)]
\item  \label{sigma-n-seq-comp} The space $\sigma_n(X)$ is sequentially compact, with a witness definable from $X$ and $n$.
\item  \label{it:prod-sigma-seq-comp} The space  $\prod_{k \in \IN} \sigma_k(X)$   is  sequentially compact, with a witness definable from $X$.
\end{enumerate}
\end{proposition}
\Pr \eqref{sigma-n-seq-comp} The proof is by induction on $n$. For $n=1$, we already know that 
$\sigma_1(X)=\hat X$ is sequentially compact with a witness definable from $X$ (use Proposition~\ref{prop:hatX-seq-comp}-\eqref{it:X-chap-2}). 
We now assume that for some integer $n\ge 1$,  
each space $\sigma_k(X)$ ($1 \le k \le n$) is  sequentially compact with a witness $\Phi_k$ definable from $X$ and $k$. Let $(F_k)_ {k \in \IN}$
be a sequence of $\sigma_{n+1}(X)$. For every $\nu \in \IN$, let 
$\mathcal A_{\nu}:= \{A \in [\IN]^{\omega} : \; \forall i \neq j \in A \; |F_i \cap F_j| =\nu\}$. 
Let $\nu_0$ be the first element of $\IN$ such that the set 
$\mathcal A_{\nu_0}$ is non-empty. One can build by induction some element  $A \in \mathcal A_{\nu_0}$,
which is definable from $X$ and $(F_n)_{n \in \IN}$. If $\nu=0$, then the subsequence 
$(F_n)_{n \in A}$ converges to $\infty$. Else, there exists $a \in X$ such that
the set $D_a := \{n \in A : a \in F_n\}$ is infinite. Build by induction some infinite subset $B$ of $A$ such that there exists   
some element $a \in X$ satisfying  $\forall n \in B \; a \in F_n$.
Let $R$ be the non-empty finite set $\cap_{n \in B}F_n$; let $p$ be the cardinal of $R$. The sequence $(F_n \backslash R)_{n \in B}$ 
lives in $\sigma_{n+1 - p}(X)$ thus, 
using the witness $\Phi_{n+1-p}$, it has an infinite  subsequence
$(F_n \backslash R)_{n \in C}$ which converges to some $L \in \sigma_{n+1 - p}(X)$. 
It follows that $(F_n \cup R)_{n \in C}$ converges to $L \cup R$ in $\sigma_{n+1}(X)$.  \\
\eqref{it:prod-sigma-seq-comp} Use Proposition~\ref{prop:hatX-seq-comp}-\eqref{polyad-den-seq} or Lemma~\ref{lem:prod-wit-seq-comp}. 
\Rp

\subsubsection{$\AC^{fin(X)}$ and closed-compactness of the space $\prod_{n \in \IN}\sigma_n(X)$}
\begin{theorem} \label{theo:pros-sigma-n}
 Let $X$ be a set. 
\begin{enumerate}[(i)]
 \item \label{it:prod-sigma-closely-comp}  $\AC^{fin(X)} \Leftrightarrow \text{``$\prod_{n \in \IN} \sigma_n(X)$   is closely compact''}$.
\item \label{it:prod-sigma-comp}  $\T_{\IN}^{fin(X)}   \Leftrightarrow \text{``$\prod_{n \in \IN} \sigma_n(X)$   is  compact''}$.
\end{enumerate}
\end{theorem}
\Pr We may assume that $X$ is infinite. In both cases, we  use Proposition~\ref{prop-im-cont} and the fact that the space  $\prod_{n \in \IN } \sigma_n(X)$ is a continuous image of 
 $\prod_{n \in \IN} \sigma_1(X)^n$, which is homeomorphic with 
 ${\hat X}^{\IN}$. \\
\eqref{it:prod-sigma-closely-comp} $\Rightarrow$: with $\AC^{fin(X)}$, ${\hat X}^{\IN}$ 
is closely compact (see Proposition~\ref{theo:acdf-comp-polyad}), 
and so is its continuous image  $\prod_{n \in \IN } \sigma_n(X)$ . 
$\Leftarrow$: if $\prod_{n \in \IN} \sigma_n(X)$  is closely compact, then so is its continuous image 
$\sigma_1(X)=\hat X$, thus $\AC^{fin(X)}$ holds.\\
\eqref{it:prod-sigma-comp} $\Rightarrow$: Using Proposition~\ref{theo:acdf-comp-polyad}-\eqref{it:hatx-puiss-alpha-1}, 
$\T_{\IN}^{fin(X)}$ implies that  ``${\hat X}^{\IN}$ is compact''. 
Using Remark~\ref{rem:sigma-n}, it follows that   ``$\prod_{n \in \IN} \sigma_n(X)$ is compact''.  $\Leftarrow$: if 
$\prod_{n \in \IN} \sigma_n(X)$   is  compact, then its closed subset ${\sigma_1(X)}^{\IN}$ is also compact, thus 
$\T_{\IN}^{fin(X)}$ holds by Proposition~\ref{theo:acdf-comp-polyad}-\eqref{it:hatx-puiss-alpha-1}.
\Rp

\section{$\AC^{fin(I)}$ and closed compactness of $B_1(I)$} \label{subsec:ac-fin2compB1}
\subsection{Dyadic representations} \label{subsec:dyadic}
\begin{notation}[{\bf binary expansion of a real number}] For every $n \in \IN$, let $\varepsilon_n := \frac{1}{2^{n+1}}$. Then  the mapping 
 $\phi : \{0,1\}^{\IN} \to [0,1]$  associating to every  $(x_n)_{n \in \IN}$ the real number 
 $\sum_n \varepsilon_n x_n$ is continuous (a  uniformly convergent series of continuous functions), onto, and 
$\phi$ has a (definable) section.
\end{notation}

\begin{theorem*}\cite[Lemma~1.1]{B-R-W} 
 Let $I$ be a set, and for every $n \in \IN$, let $I_n:= \{n\}  \times I$.  Let $F$ be a closed subset of $[0,1]^I$. 
Consider  the power mapping $g:=\phi^I :  \{0,1\}^{\IN \times I} \to [0,1]^I$. 
For every $n \in \IN$, let $j_n :  \{0,1\}^{I_n} \to  \{0,1\}^{\IN \times I}$ be the canonical inclusion mapping.
Let $Z:= g^{-1}[F]$ and, for every $n \in \IN$, let $Z_n := j_n^{-1}[Z]$: thus $Z$ is a closed subset of $\prod_{n \in \IN}Z_n$ 
and $g: Z \to F$ is continuous, onto, with a definable section 
\begin{enumerate}[(i)]
\item \label{it:dyadic2} If $F \subseteq B^+_1(I)$, then for every $n_0 \in \IN$, $Z_{n_0} \subseteq \sigma_{M_{n_0}}(I_{n_0})$ where 
$M_{n_0}:= \lfloor \frac{1}{\varepsilon_n} \rfloor$ (the integral part of $\frac{1}{\varepsilon_{n_0}}$), thus $F$ is the continuous image of some closed subset of ${\hat I}^{\IN}$. 
\item \label{it:dyadic3} If $F \subseteq \ell^0(I)$, then for every $n_0 \in \IN$, $Z_{n_0} \subseteq \{0,1\}^{(I_{n_0})}$
\end{enumerate}
 \end{theorem*}
\Pr Let $n_0 \in \IN$ and let  $(x^{n_0}_i)_{i \in I} \in Z_{n_0}$;  let $x=(x^n_i)_{n \in \IN, i \in I} \in Z$
such that  $j_{n_0}((x^{n_0}_i)_{i \in I})=x$.
\eqref{it:dyadic2}  Since $F \subseteq B^+_1(I)$,  $\sum_{i,n} \varepsilon_n x^n_{i} \le 1$, thus
$\sum_{i \in I} \varepsilon_{n_0} x^{n_0}_{i} \le 1$; it follows that the set $\{i \in I :  x^{n_0}_{i}=1 \}$ has a cardinal 
$\le  \frac{1}{\varepsilon_{n_0}}$.\\
\eqref{it:dyadic3} Since  $F \subseteq \ell^0(I)$,   $(x^{n_0}_i)_{i \in I} \in \ell^0(I) \cap \{0,1\}^I$
thus  $Z_{n_0} \subseteq \{0,1\}^{(I_{n_0})}$. 
\Rp

\begin{remark} \label{rem:B+&B} The mapping $- : B^+_1(I) \times B^+_1(I) \to B_1(I)$ is  continuous and onto thus
$B_1(I)$ is also the continuous image of a closed subset of ${\hat I}^{\IN}$.
\end{remark}

\begin{remark} Aviles  (\cite{Avi07}) proved that $B^+_1(I)$ -and thus $B_1(I)$-  
is a continuous image of  ${\hat I}^{\IN}$ (and not only of a closed subset of ${\hat I}^{\IN}$).
\end{remark}

\subsection{Another equivalent of $\AC^{fin(I)}$} \label{subsec:acfin2-u-eber}
\begin{theorem} \label{theo:ac-fin2-ball-hilb-base} Let $I$ be a set.
 \begin{enumerate}[(i)]
 \item \label{it:acf&B1} $\AC^{fin(I)} \Leftrightarrow \text {``$B_1(I)$  is closely compact.''}$    
Moreover,  a witness of closed compactness on $B_1(I)$ is 
definable from $I$ and a choice function for non-empty finite subsets of $I$ and conversely. 
\item   \label{it:tych&B1} $\T_{\IN}^{fin(I)}$ implies that  $B_1(I)$ is compact. 
\item  \label{it:ZF&B1} The space  $B_1(I)$ is sequentially compact, with a witness definable from $I$. 
\end{enumerate}
\end{theorem}
\Pr Using the previous Theorem and Remark~\ref{rem:B+&B}, consider some sequence $(M_n)_{n \in \IN}$,  some closed subset
$Z$ of $\prod_{n \in \IN} \sigma_{M_n}(I)$ and some continuous onto mapping $g : Z \to B_1(I)$, with  a (definable) section. \\
\eqref{it:acf&B1}  $\Rightarrow$: Using $\AC^{fin(I)}$ and Theorem~\ref{theo:pros-sigma-n}-\eqref{it:prod-sigma-closely-comp}, $Z$  is  closely compact
 thus $B_1(I)=g[F]$ is closely compact.  $\Leftarrow$: If $B_1(I)$ is closely compact, then $\hat I$ (which is 
 a closed subset of $B_1(I)$ -see Example~\ref{ex:hatX_unif-eb} in Section~\ref{subsec:alex}-) is also closely compact. \\
  \eqref{it:tych&B1}  Using  $\T_{\IN}^{fin(I)}$ and Theorem~\ref{theo:pros-sigma-n}-\eqref{it:prod-sigma-comp},  $Z$ is  compact thus $B_1(I)=g[F]$ is also compact.  \\
\eqref{it:ZF&B1} By Proposition~\ref{prop:seq-comp}-\eqref{it:prod-sigma-seq-comp}, 
$Z$ is sequentially compact and $g$ is continuous with  a section, thus Lemma~\ref{lem:imag-cont-seq-comp}
implies that $B_1(I)$
is also sequentially compact with a witness definable from $I$. 
\Rp

\begin{corollary} \label{cor:equ-B1-comp} Given a set $I$, the following statements are equivalent:
 \begin{enumerate}[(i)]
\item \label{it:tych-omega-fin-1} $\Uwo_{\IN}^{fin(I)}$
\item \label{it:tych-omega-fin-2} $\T_{\IN}^{fin(I)}$
 \item \label{it:tych-omega-fin-3} The space  $B_1(I)$ is    compact.
\item \label{it:tych-omega-fin-4} For every sequence $(F_n)_{n \in \IN}$ of finite subsets of $I$, the space  $B_1(\cup_{n \in \IN} F_n)$ is    compact.
\end{enumerate}
\end{corollary}
\Pr \eqref{it:tych-omega-fin-1}   $\Rightarrow$ \eqref{it:tych-omega-fin-2} is easy and 
\eqref{it:tych-omega-fin-2}   $\Rightarrow$ \eqref{it:tych-omega-fin-3} follows from Theorem~\ref{theo:ac-fin2-ball-hilb-base}. 
\eqref{it:tych-omega-fin-3}   $\Rightarrow$ \eqref{it:tych-omega-fin-4} is easy. We show that 
\eqref{it:tych-omega-fin-4}   $\Rightarrow$ \eqref{it:tych-omega-fin-1}. 
 The idea of the implication is in \cite[th.~9 p.~16]{Fo-Mo}: we sketch it for sake of completeness. 
Let $(F_n)_{n \in \IN}$ be a disjoint sequence of non-empty finite sets
of $I$. Let us show that $D:=\cup_{n \in \IN}F_n$ is countable.
The Hilbert spaces $H:=\ell^2(D)$ and $\oplus_{\ell^2(\IN)} \ell^2(F_n)$ are isometrically isomorph. 
For every $n \in \IN$, let $\varepsilon_n : |F_n| \to ]0,1[$ be 
a strictly increasing mapping such that 
$\sum_{n \in \IN} \sum_{0 \le i < |F_n|}\varepsilon_n(i)^2=1$. 
For every $n \in \IN$, let 
$\tilde F_n :=\{x \in B_{H} : \; \forall m \neq n, x_{\restriction F_m}=0 \text{ and } x_{\restriction F_n} \text{ is one-to-one from } F_n \text{ onto } rg(\varepsilon_n)\}$.
Then  each $\tilde F_n$ is a weakly closed subset of $B_{\ell^2(D)}$ and 
the sequence $(\tilde F_n)_{n \in \IN}$ satisfies the finite intersection property. 
The  compactness  of $B_2(D)$ implies that 
$Z:=\cap_{n \in \IN} \tilde F_n$ is non-empty. Given an element $f=(f_n)_{n \in \IN}$ of $Z$,  each $f_n$ defines a well-order on the finite set $F_n$, thus $\cup_{n \in \IN} F_n$ is countable. 
\Rp

\begin{remark} For $I=\mathcal P(\IR)$, none of the equivalent statements in Corollary~\ref{cor:equ-B1-comp} 
is provable in $\ZF$.
Indeed, there is a model of $\ZF$ where there exists a sequence $(P_n)_{n \in \IN}$ of pairs of subsets of $\IR$
such that $\prod_{n \in \IN} P_n$ is empty.
\end{remark}

Thus, the  statement {\em ``The closed unit ball of $\ell^2(\mathcal P(\IR))$ is weakly compact.''} is not provable in $\ZF$. 

\begin{corollary} \label{cor:acdf&Ahb} The following statements are equivalent: $\AHb$, $ \AHbf$, $\ACDF$.
\end{corollary}

\subsection{Consequences}
\begin{corollary} \label{cor:B1deLO}
If a set $I$ is linearly orderable, then $B_1(I)$ is compact.
\end{corollary}
\Pr If $I$ is linearly orderable, then $\AC^{fin(I)}$ holds. 
\Rp

For every ordinal $\alpha$, the set $\mathcal P(\alpha)$ is linearly orderable, thus  $\AC^{fin(\mathcal P(\alpha))}$ holds. 
In particular, $\IR$ is equipotent with $\mathcal P(\IN)$ so the closed unit ball of $\ell^2(\IR)$ is closely compact.  
This solves Question~3 of \cite{Mo07}.

\begin{question}
Does  $\AC^{fin}$ imply $\AH$?  What is the power of the statement {\em ``Every Hilbert space
has a hilbertian basis''}? Is  this statement provable in $\ZF$? Does it imply $\AC$?
\end{question}

\section{$\ACD$ and Eberlein spaces} \label{sec:dc-and-comp}
Given a set $I$, a closed subset $F$ of $[0,1]^I$ is {\em $I$-Corson} if every element $x \in F$ has a countable support. 
\subsection{Sequential compactness of $I$-Eberlein spaces}
Given a set $I$, denote by $count(I)$ the set of finite or countable subsets of $I$.
Consider the  following consequence of $\ACD$: \\
\noindent $\Uwo_{\IN}^{count(I)}$: {\em ``Every countable union of countable subsets of $I$ is countable.''}

\begin{proposition} \label{prop:seq-comp-corson} Let $I$ be a set and let $F$ be a closed subset of $[0,1]^I$. 
\begin{enumerate}[(i)]
\item  \label{it:strong-eber-comp2} $\Uwo^{fin(I)}_{\IN}$ implies that:
\subitem  -if $F \subseteq [0,1]^{(I)}$, then $F$ is sequentially compact.
\subitem  -if  $F \subseteq c_0(I)$, then $F$ is $I$-Corson.
\item $\Uwo^{count(I)}_{\IN}$ implies that if  $F$ is $I$-Corson,  
then $F$ is    sequentially compact. Thus $\ACD$ implies that every Eberlein space is sequentially compact.
\item $\AC^{fin(I)}$ implies that if $F \subseteq [0,1]^{(I)}$,  
then $F$ is sequentially compact and has a witness of sequential 
compactness.
\end{enumerate}
\end{proposition}
\Pr \eqref{it:strong-eber-comp2} If $F \subseteq [0,1]^{(I)}$, then $F$ is sequentially compact using 
the fact that $[0,1]^{\IN}$ is sequentially compact (see Example~\ref{ex:seq-compact}).
\Rp 

\subsection{Countable product of finitely restricted spaces}
\begin{theorem} \label{theo:acd-stron-eber-comp}
 Let $F$ be a closed subset of $[0,1]^I$ which is contained in $[0,1]^{(I)}$. Then
$\ACD$ implies that $F$ is compact. In particular, $\ACD$ implies that every strong Eberlein  space is compact.
\end{theorem}
\Pr  For every $n \in \IN$, recall that the subset  
$\sigma_n(I) :=\{x \in [0,1]^{(I)} : \; |supp(x)| \le n\}$ is compact (see Proposition~\ref{prop:sigma-n-comp}). 
Let  $\mathcal F$ be a filter of the lattice of 
closed subsets of  $F$.
If there exists an integer $n$ such that  $\sigma_n(I) $ is $\mathcal F$-stationar, 
then $\cap \mathcal F$ is non-empty by compactness of $\sigma_n(I)$ and using Remark~\ref{rem:stat-comp}. 
Else, using $\ACD$, consider a sequence  $(F_n)_{n \in \IN}$ of closed subsets of $F$ belonging to  
$\mathcal F$ such that for every $n \in \IN$, $F_n \cap \sigma_n(I) =\varnothing$.
Re-using  $\ACD$, choose for every  $n \in \IN$ an element 
$x_n \in F_n$. 
A new use of $\ACD$ and Proposition~\ref{prop:seq-comp-corson}-\eqref{it:strong-eber-comp2} implies  the existence of  some infinite subset $A$ of $\IN$ such that 
  $(x_n)_{n \in A}$ converges to some element $x \in F$.
Then, for every  $n \in \IN$, $x \in F_n$ (which is disjoint with $\sigma_n(I)$) 
so the element $x$ of $F$ has an infinite support: this is contradictory! 
\Rp

\begin{remark}
 It does not seem provable in $\ZF$ that  every closed subset of $[0,1]^I$ contained in $[0,1]^{(I)}$ is compact, or has a witness of sequential compactness.
\end{remark}

\subsection{Countable products of strong Eberlein spaces}
\begin{theorem} \label{theo:prod-strong-eber}
Let  $(I_n)_{n \in \IN}$ be a sequence of pairwise disjoint sets and let 
$I:=\cup_{n \in \IN} I_n$.
For every  $n \in \IN$,  let $F_n$ be a closed subset of  $\{0,1\}^{(I_n)}$. 
Let  $F$ be the closed subset $\prod_{n \in \IN} F_n$ of  $\{0,1\}^I$.
Then $\ACD$ implies that $F$ is   sequentially compact and compact.
\end{theorem}
\Pr  Using $\ACD$, $F$ is sequentially compact: given a sequence $(x_n)_{n \in \IN}$ of $F$, $\ACD$ implies that for every $n \in \IN$, the support $D_n$ of $x_n$ is countable, thus re-using $\ACD$, the set $D:=\cup_{n \in \IN}D_n$ is also countable; since each $x_n$ belongs to 
$[0,1]^D \times \{0\}^{I \backslash D}$, and since $[0,1]^D$ is sequentially compact (see Example~\ref{ex:seq-compact}), it follows that $(x_n)_{n \in \IN}$ has an infinite 
 subsequence which converges in $F$. Using $\ACD$ and Theorem~\ref{theo:acd-stron-eber-comp}, each $F_n$ is  compact. Using $\ACD$ 
and Remark~\ref{rem:acd-seq-comp}-\eqref{it:acd-prod-seq-comp}, it follows that $F$ is compact.   
\Rp 

\subsection{$\ACD$ and  Eberlein closed subsets of  $[0,1]^I$} \label{subsec:acd2eberl}
\begin{corollary} \label{cor:acd2comp-eberlein}
 $\ACD$ implies that every Eberlein space is  both sequentially compact and compact. 
\end{corollary}
\Pr Let $X$ be an  Eberlein space. Then $X$  is sequentially compact by Proposition~\ref{prop:seq-comp-corson}-\eqref{it:strong-eber-comp2}.
Let $I$ be a set such that $X$ is homeomorphic with a closed subset $F$ of $[0,1]^I$, with $F \subseteq [0,1]^{(I)}$.
Using the Theorem of Section~\ref{subsec:dyadic},   there exists a family $(Z_k)_{ k \in \IN}$ such that for each $ k \in \IN$, 
$Z_k$ is $I$-strong Eberlein, and such that $F$ is the continuous image of a closed subset $Z$ of $\prod_{k \in \IN} Z_k$. 
Using $\ACD$ and Theorem~\ref{theo:prod-strong-eber}, the space $\prod_{k \in \IN} Z_k$ is sequentially compact. 
With $\ACD$,  each $Z_k$ is compact;  with Remark~\ref{rem:acd-seq-comp}-\eqref{it:acd-prod-seq-comp} and $\ACD$, it follows that $prod_{k \in \IN} Z_k$ (and thus its continuous image $F$) is also compact. 
\Rp

Recall that (see Theorem~\ref{theo:ac-fin2-ball-hilb-base}-\eqref{it:ZF&B1}) every $I$-uniform Eberlein space
is sequentially compact, with a witness definable from $I$.

\begin{question} Does $\ACD$ implies that every Eberlein space is  sequentially compact with a witness?
\end{question}

\subsection{Convex-compactness and the Hahn-Banach property}
Given  a set $I$, say that a subset $F$ of $\IR^I$ is {\em convex-compact} if for every set $\mathcal C$ of closed convex 
subsets of $\IR^I$ such that $\{F \cap C : \; C \in \mathcal C\}$ satisfies the {\em FIP}, $C \cap \bigcap \mathcal C$
is non-empty; moreover, if there is a mapping associating to each such  $\mathcal C$ an element of $C \cap \bigcap \mathcal C$,
then say that $F$ is {\em closely} convex-compact. 
Given a topological vector space  $E$, say that $E$ satisfies the {\em continuous Hahn-Banach property} if, for every continuous sublinear
functional $p : E \to \IR$, for every vector subspace $F$ of $E$, and every linear functional
$f : F \to  R$ such that $f \le  p|F$ , there exists a linear functional $g : E\to \IR$ that
extends $f$ and such that $g \le  p$. Moreover, if there is a mapping associating to each such $f$ some $g$ satisfying the
previous conditions, then say that $E$ satisfies the {\em effective} continuous Hahn-Banach property.

\begin{theorem*}[\cite{Fo-Mo}, \cite{Do-Mo}] Given a set $I$, the normed space $\ell^0(I)$ satisfies the effective continuous Hahn-Banach property.
For every real number $p \in [1,+\infty[$, $\ell^p(I)$ satisfies the effective continuous Hahn-Banach property.
\end{theorem*}

\begin{corollary*} For every set $I$, and every real number $p \in [1,+\infty[$, $B_p(I)$ is closely
convex-compact.
\end{corollary*}
\Pr The continuous dual of $\ell^0(I)$ is (isometrically isomorphic with) $\ell^1(I)$ and, for every 
$p \in ]1,+\infty[$,  the continuous dual of $\ell^p(I)$ is $\ell^q(I)$ where $\frac{1}{p} + \frac{1}{q}=1$ thus
$1<q<+\infty$.
We end the proof using  the fact (see \cite{Fo-Mo})  that if a normed space $E$ satisfies the (effective) continuous Hahn-Banach property, then the closed unit ball
of the continuous dual $E'$ is ({\em closely}) weak* compact.
\Rp

\begin{question} Given a set $I$, and a closed convex subset $C$ of $[0,1]^I$ which is $I$-Eberlein,
is $C$ convex-compact (in $\ZF$)?
\end{question}

\begin{question} Same question if the closed convex subset $C$ of $[0,1]^I$ is $I$-Corson.
\end{question}

\bibliographystyle{abbrv}
\bibliography{../../biblio}
\end{document}